\documentclass[reqno]{amsart}
\usepackage{graphicx}
\usepackage{indentfirst,csquotes}

\usepackage{amssymb,amsthm,amsmath}
\usepackage{xcolor,paralist,hyperref,titlesec,fancyhdr,etoolbox}

\usepackage{amsthm}
\theoremstyle{remark}
\newtheorem{remark}{Remark}

\titleformat{\section}
  {\normalfont\Large\bfseries\centering}
  {\thesection}{1em}{}

\titlespacing*{\subsection}
{0pt}      % left indent (keep 0pt => left aligned)
{1.0ex}    % space before
{0.8ex}    % space after (increase this to separate from paragraph)

\titlespacing*{\subsubsection}
{0pt}
{0.8ex} 
{0.6ex}
\titleformat{\subsection}
  {\normalfont\large\bfseries}  % style (left aligned by default)
  {\thesubsection}{1em}{}

\hypersetup{ colorlinks=true, linkcolor=black, filecolor=black, urlcolor=black }

\usepackage{lipsum}

\begin{document}
\title{LPV Updates for Sequentially Linearized Moving Horizon Estimation of Nonlinear Systems} %%%%%%%%%%%%
\author[Initial Surname]{JIAXIN JI, JAN HEILAND, Dimitrios S. Karachalios and Hossam S. Abbas}
%\address{Address}
%\email{example@mail.com}

% \let\thefootnote\relax
% \footnotetext{MSC2020: Primary 00A05, Secondary 00A66.} %%%%%%%%%%

\begin{abstract}
 \emph{Moving horizon estimation} (MHE) provides high precision state estimation for nonlinear systems, but it is often limited by the substantial computational demands of solving a nonlinear optimization problem at every sampling step. To address this issue, we develop an efficient MHE scheme based on \emph{linear parameter-varying} (LPV) formulation, where the scheduling parameters are given by the estimated states of the system and used to construct inexact Jacobians.
 %for the \emph{quadratic programming} (QP) subproblem.
 %update inexact Jacobians used in the \emph{quadratic programming} (QP) subproblem. 
Due to the LPV representation, the Jacobian can be pre-specified offline in a structured form and then updated in the \emph{quadratic programming} (QP) subproblem, which reduces computational cost commonly used in standard \emph{nonlinear programming} (NLP) systems.
 We illustrate the performance by numerical simulations.
\end{abstract} %%%%%%%%%
\maketitle
% \bigskip

% \noindent \lipsum[1] \cite{1}

% $\,$

% $\,$
\section{Introduction}
For many control applications, system behavior cannot be characterized only by inputs and outputs. Instead, internal state variables are required to capture the dynamics. Since these states are typically not fully measurable and sensor data are affected by noise, state estimation is required. 
A widely used estimation method for nonlinear systems is Moving Horizon Estimation (MHE) \cite{VBH01}, which can achieve high accuracy and incorporate constraints by solving an optimization problem over a sliding window of recent measurements. However, MHE requires solving a nonlinear program within the sampling instant, which may lead to substantial computational cost and limits its use in real-time settings. For high-dimensional or strongly nonlinear systems, solving such nonlinear optimization problems online becomes even more complex. % To address this issue, real-time iteration (RTI) approaches have been proposed, such as the RTI scheme in ACADO developed by Diehl et al.~\cite{KD11}. However, the Jacobian evaluations still need to be performed for the nonlinear optimization problem.
To address this issue, a variety of linearization-based MHE methods have been developed. One representative approach is the real-time iteration (RTI) scheme developed by Diehl et al. and implemented within the ACADO framework \cite{KD11}. By performing only one generalized Gauss-Newton iteration at each sampling instant, the RTI scheme substantially reduces the online computational burden. Related finite-iteration approaches for discrete-time systems were studied by Alessandri and Gaggero \cite{AleG17}, who considered single and multi iteration descent schemes based on gradient, conjugate gradient and Newton directions. However, reducing the number of descent iterations does not eliminate the cost of evaluating the derivatives of the nonlinear MHE problem. Further reductions have been obtained by reusing derivative information. In particular, the zero-order MHE method of Baumgärtner et al. \cite{BaFH21} combines a Gauss-Newton formulation with fixed sensitivity approximations and reused matrix factorizations, thereby reducing repeated online sensitivity evaluations for large-scale nonlinear processes. More recently,  Schiller and Müller \cite{SchM22} showed that, under an appropriate estimator design, useful estimation performance and robust stability properties can be maintained with incomplete optimization. These methods all construct the optimization updates from local approximations of the original nonlinear optimization problem.

In this work, we propose an MHE scheme based on Linear Parameter-Varying (LPV)  by considering the recent Sequential Quadratic Programming (SQP)-LPV embedding approach in \cite{KA24} for solving the model predictive control (MPC) problem, with the goal of reducing the computational burden of nonlinear MHE. A quasi-LPV representation preserves a linear state space structure while describing the nonlinear dependence through scheduling variables, which are
functions of internal system signals. Systematic methods for embedding control-affine nonlinear systems into LPV representations were developed, for example, by Abbas et al. \cite{AbbT14}. MHE for LPV and quasi-LPV systems has been investigated previously from different perspectives. Alessandri et al. \cite{AlZZ20} proposed optimistic and pessimistic MHE formulations for quasi-LPV systems with uncertain scheduling variables, where the scheduling trajectory is jointly estimated with the state trajectory. Robust stability guarantees for quasi-LPV MHE were subsequently established in \cite{ArAZ23}. Miyakoshi et al. \cite{MiKK24} proposed an LPV-based MHE framework for jointly estimating vehicle states and cornering stiffness. By fixing either the state variables or the varying parameters, the original coupled estimation problem is decomposed into a sequence of convex optimization problems. However, exploiting the LPV structure directly in the numerical solution of the nonlinear MHE optimization problem has received little attention. In particular, to the best of our knowledge, existing LPV-based MHE approaches do not utilize the scheduling-dependent LPV structure to construct structured optimization updates for solving nonlinear MHE. The proposed work exploits the LPV representation at the optimization level. Specifically, we use an LPV representation to derive the structured Jacobian directly, and then update the values by scheduling parameters in LPV.

% Based on this LPV representation, Karush–Kuhn–Tucker (KKT) conditions associated with the Nonlinear Programming (NLP) of the MHE optimization problem are approximated and formulated in terms of the scheduling parameter and solved using an inexact Newton-type update, leading to a Sequential Quadratic Programming (SQP) framework with a structured QP subproblem at each iteration, which can be systematically updated via the scheduling parameter, eliminating the need to compute gradients or Hessians in each iteration of conventional SQP algorithms for NLP. The resulting SQP-MHE scheme preserves the estimation performance of nonlinear MHE while substantially reducing computation through linear constraints in the LPV model and the quadratic cost formulation.
In section ~\ref{sec:foundation}, we first introduce the nonlinear system model and the MHE formulation used in this work and then present the quasi-LPV reformulation. 
In section ~\ref{sec:new scheme}, we derive the proposed LPV-SQP-MHE framework. The nonlinear MHE problem is first expressed in an LPV form, in which the nonlinear dynamics is expressed as a parameter-dependent linear structure. Then the Karush–Kuhn–Tucker (KKT) conditions \cite{Ber97} associated with the MHE optimization problem are approximated, formulated in terms of the scheduling parameter. %, and solved,
The resulting KKT conditions are solved within an SQP framework, where a structured Quadratic Programming (QP) subproblem is obtained at each iteration. %leading to a Sequential Quadratic Programming (SQP) framework with a structured QP subproblem at each iteration. 
In particular, the systematic update through the LPV parameter partially circumvents the computation of Jacobians as in conventional SQP schemes. Specifically, we introduce two algorithms based on different linearization strategies to provide the QP subproblem. 
Section ~\ref{sec:simulation} reports simulation results on a two-link arm robot manipulator, whose model is based on the formulation in \cite{JiJM25}. We assess the accuracy and runtime of the estimation of two proposed approaches. Moreover, we compare the performance of two LPV-SQP-MHE approaches with that of the standard nonlinear MHE solved by CasADi/IPOPT \cite{AnGHRD18}. 
Finally, Section~\ref{sec:conclusion} concludes the paper and outlines directions for future work.

Throughout the paper, we use the following notation. 
$ I$ and $0$ denote the identity matrix and the zero matrix with compatible dimensions, respectively. 
For a symmetric matrix $W\succ 0$, we use the weighted norm $\|x\|_W^2 := x^\top W x$.
Furthermore, for a matrix $X$, we use  $\|X\|_F$ to denote the Frobenius norm.
The Kronecker product is written as $\otimes$.

\section{Problem formulation}
\label{sec:foundation}
%--------------------------------------nonlinear systme 

In this work, we consider a nonlinear system defined by the state-space representation:
\begin{equation}
\begin{aligned}
 \label{eq:nonlinear}
x_{k+1} &= f(x_k, u_k)+w_k, \\
\quad y_k &= h(x_k)+v_k,
\end{aligned}
\end{equation}
where \(k\in\mathbb N\) denotes the discrete-time sampling index, $x_k \in \mathbb{R}^{n}$ is the state variable with \eqref{eq:nonlinear},
while $u_k \in \mathbb{R}^{m}$ is the known input, and $y_k \in \mathbb{R}^{n_y}$ is the output of the system. The functions $f : \mathbb{R}^n \times \mathbb{R}^m \to \mathbb{R}^{n}$ and $h : \mathbb{R}^n \to \mathbb{R}^{n_y}$ are vector-valued functions. $w_k$ and $v_k$ denote the process noise and measurement noise, respectively.

% The process noise $w_k$ and the measurement noise $v_k$ are constrained within the set $\mathcal{W}$ and $\mathcal{V}$, defined as 
% \begin{equation}
% \label{eq:general noise model}
% \begin{aligned}
%     \mathcal{W}:= \{ w \in \mathbb{R}^n \quad \text{s.t} \quad \| w\|_0 \leq \theta_w \} \\
%      \mathcal{V}:= \{ v \in \mathbb{R}^{n_y} \quad \text{s.t} \quad \| v\|_0 \leq \theta_v \}
% \end{aligned}
% \end{equation}
% where $\theta_w >0 $ and $\theta_v >0 $ are known scalars. 
%-------------------------------------- LPV embedding  here or later?

% Embedding of the NL system in \eqref{eq:nonlinear} in an LPV representation corresponds to constructing an equivalent system of the following
% \begin{equation}
% \begin{aligned}
% \label{LPV}
% x_{k+1} &= A(\rho_k)x_k + B(\rho_k)u_0 +w_k, \\
% y_k &= C(\rho_k)x_k +v_k, 
% \end{aligned}
% \end{equation}
% where $\rho_k \in  \mathbb{R}^{n_p}$ is the scheduling variable, and there exists a function $\eta$ called scheduling map, such that $ \eta(x_k, u_k) = \rho_k$. For now, we let the scheduling variable $\rho_k=x_k$, which assumes that \eqref{eq:nonlinear} is input-affine.

%--------------------------------------standard MHE
\subsection{Moving Horizon Estimation}
\begin{figure}[t]
  \centering
  \includegraphics[width=0.8\textwidth, height=0.5\textwidth]{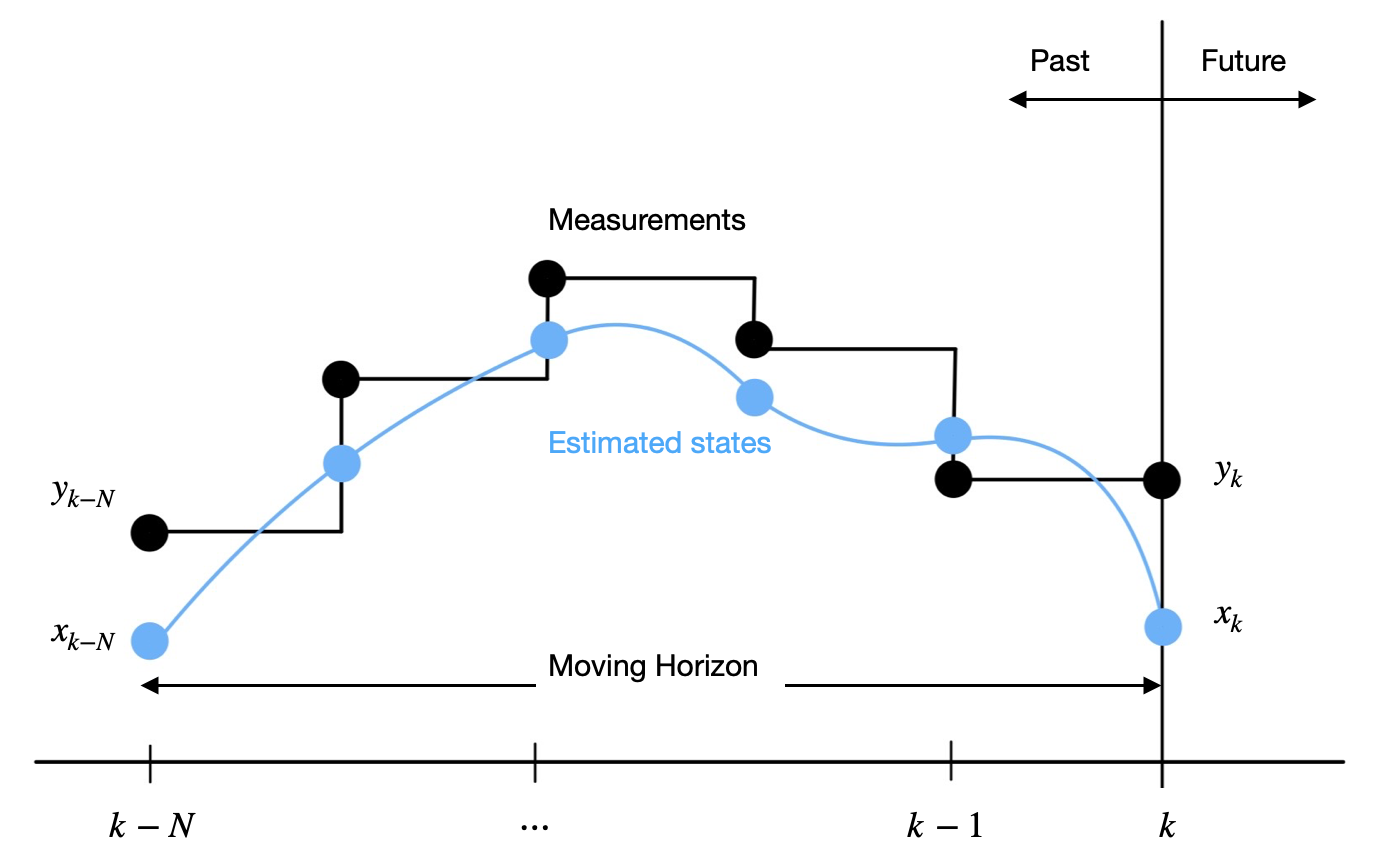} 
  \caption{MHE framework}
  \label{fig:mhe}
\end{figure}

In this subsection, the MHE problem is first written as a nonlinear program, where the decision vector contains the initial state and noise estimates over the moving horizon. Let $N$ denote the length of an estimation window. 
Consider data within the finite time window that move over time, and MHE uses only the measurements and control inputs within the window in each estimation step to estimate the current state of the system. The horizon principle of MHE is illustrated in Fig.\ref{fig:mhe}.

%---------figure 

At each sampling time $k$, the information consists of the $N$ historical data. For $i=0,\ldots,N$, let $\chi_{k-N+i|k}$ denote the estimation of the state in the time step $k-N+i$ based on the available measurements up to the time step $k$, and let $y_{k-N+i}$ denote the measurement at time step $k-N+i$. For notational simplicity, we define 
\begin{equation*}
\chi_i\equiv \chi_{k-N+i|k}, \quad
    y_i := y_{k-N+i}.
\end{equation*}
In particular, the initial state of the current estimation window is $\chi_0\equiv \chi_{k-N|k}$.
% For simplicity, we index variables locally within the window by $i=0,\cdots, N$. In particular, the window initial state $\chi_0$ denotes the state at the beginning of the current $k$-th window, i.e. $\chi_0\equiv \chi_{k-N|k}$. Also, the local window variables $\chi_i$ and $y_i$ correspond to the global-time variables:
% \begin{equation*}
%     \chi_i\equiv \chi_{k-N+i|k}, \quad y_i \equiv  y_{k-N+i|k} 
% \end{equation*}
The MHE problem is defined as:
\begin{equation}
\label{eq:MHE standard}
\begin{aligned}
\min_{\chi_0,\omega,\nu}f(z)= \min_{\chi_0,\omega,\nu}
\frac12||\chi_0-\bar x_0||_{P}^2 +
\frac12\sum_{i=0}^{N-1}||\omega_{i}||_{Q}^2 +
\frac12\sum_{i=0}^{N}||\nu_{i}||_{R}^2 \
\end{aligned} ,
\end{equation}
\vspace{-0.2 cm}
    subject to
\vspace{-0.1 cm}
    \begin{equation}
    \label{eq:standard constaints}
    \begin{aligned}
\chi_{i+1} & =f(\chi_i,u_i) + \omega_{i},\quad i=0,\dots,N-1 ,\\ %A(\rho(x_{k}))x_{k} + B(\rho(x_{i}))u_{i} 
 y_i &= h(\chi_i)+ \nu_{i},\quad i=0,\dots,N ,    %
    \end{aligned}
    \end{equation}
where $\bar{x}_0$ is the prior information on the initial state over one moving horizon, the arrival cost is weighted by a positive weighting matrix $P$, and the symmetric positive definite matrices $Q$ and $R$ represent the weighting matrices of stage cost. At the beginning of the experiment (small $k$), only a limited number of measurements is available. %Therefore, we use a shorter window initially and increase the window length until enough  have been collected. Afterwards, the horizon length N can be kept fixed.
Therefore, we start the state estimation once the horizon length is reached, i.e. at $k=N$. From this time step, the estimation is performed at every time step $k$, while the horizon length is kept fixed.

%--------------------------------------LPV representation
\subsection{Quasi-LPV representation}
A LPV representation describes a nonlinear system by a linear structure depending on time-varying scheduling variables. If scheduling variables depend on internal signals, such as state, the model is referred to as a quasi-LPV representation; see \cite{JiJM25}.
% Embedding of the NL system in \eqref{eq:nonlinear} into an LPV representation corresponds to constructing an equivalent system of the following form:
We assume that the nonlinear system in \eqref{eq:nonlinear} admits an LPV representation i.e., there exist a well-defined scheduling map $\eta(x_k, u_k) = \rho_k$ and matrix-valued functions $A(\cdot)$, $B(\cdot)$, and $C(\cdot)$ such that
\begin{equation}
\begin{aligned}
\label{eq:LPV}
x_{k+1} &= A(\rho_k)x_k + B(\rho_k)u_k +w_k, \\
y_k &= C(\rho_k)x_k +v_k, 
\end{aligned}
\end{equation}
where $\rho_k \in  \mathbb{R}^{n_p}$ is the scheduling variable. For now, we let the scheduling variable $\rho_k=\rho_k(x_k)$, which assumes that \eqref{eq:nonlinear} is input-affine. If we drop the dependency (e.g. by freezing the parameters at the values $\bar \rho_k$ computed from the estimated state trajectories obtained at the previous time step of the MHE problem), then the model becomes linear in $(x_k,u_k)$.
Using the LPV representation in \eqref{eq:LPV}, we reformulate the nonlinear MHE problem by replacing the original dynamics with the scheduling-dependent linear structure. Specifically, we keep the same objectives in \eqref{eq:MHE standard} as in nonlinear form and express the constraints in \eqref{eq:standard constaints} in quasi-LPV form:
\begin{equation}
\label{eq:nonlinear lpv constrains}
\begin{aligned}
\chi_{i+1} & =A(\rho(\chi_i))\chi_i + B(\rho(\chi_i))u_i+ \omega_{i},\quad i=0,\dots,N-1 , \\ %A(\rho(x_{k}))x_{k} + B(\rho(x_{i}))u_{i} 
y_i &= C(\rho(\chi_i))\chi_i+ \nu_{i},\quad i=0,\dots,N .     %
\end{aligned}
\end{equation}
% where the scheduling variable is given by the map $\rho_i=x_i$. For a fixed scheduling trajectory $\rho^\ell$, the dynamic becomes linear in $(x,u)$.
%--------------------------------------------------------------------------------------------------------------------------------------------------------------------------LPV-SQP-MHE scheme
\section{LPV-SQP-MHE scheme}
\label{sec:new scheme}
First, we group all decision variables over the estimation horizon into the
\begin{equation*}
    \bar{x}=
\begin{bmatrix}
    \chi_0 \\ \chi_1 \\ \vdots \\ \chi_N
\end{bmatrix}, 
    \bar{\omega}=
\begin{bmatrix}
    \omega_{0} \\ \omega_{1} \\ \vdots \\ \omega_{N-1}
\end{bmatrix},
    \bar{\nu}=
\begin{bmatrix}
    \nu_{0} \\ \nu_{1} \\ \vdots \\ \nu_{N}
\end{bmatrix},
%     \bar{\rho}=
% \begin{bmatrix}
%     \rho_1 \\ \rho_2 \\ \vdots \\ \rho_N
% \end{bmatrix}
\end{equation*}
For $z=(\bar x,\bar \omega, \bar \nu)$ as decision variable, we can write the overall constrained optimization
problem as
\begin{equation}
  \min_z f(z), \quad \text{s.t. } h(z)= 0,
\end{equation}
where
\begin{equation}
  h(z) = \begin{pmatrix}
\mathcal A(\rho( z )) z
-\mathcal B(\rho( z))u  \\
y - \mathcal C(\rho(z)) z
\end{pmatrix} ,
\end{equation}
represents the discrete system, where matrices $\mathcal{A}$, $\mathcal{B}$ and $\mathcal{C}$ are given in Appendix.
%------------------Lagrangian
Then, the Lagrangian is defined as
% \begin{equation}
%   \mathcal L(x,\lambda) = f(x) + \lambda^T h(x)
% \end{equation}
\begin{equation}
  \mathcal L(z,\lambda) = f(z) + \lambda^\top h(z) 
\end{equation}
%and the "$h$"-part of $\nabla_x \mathcal L(x, \lambda)$ is obtained as
and the "$h$"-part of $\nabla_z \mathcal L(z, \lambda)$ is obtained as
\begin{equation}\label{eq:gradientx-state-eqs}
\begin{split}
  \nabla_z (\lambda^\top h(z)) &= \begin{pmatrix}
   \lambda^\top \nabla_z\bigl (\mathcal A(\rho(z))\,z - \mathcal
  B(\rho( z))\,u \bigr) \\
   \lambda^\top \nabla_z\bigl (      y-\mathcal{C}(\rho(z)) z \bigr)
  \end{pmatrix} \\
 &= \begin{pmatrix}
      \lambda^\top \bigl (\mathcal A(\rho( z))\,[\cdot] +
                           (\nabla_\rho \mathcal A(\rho(z))\nabla_z
                         \rho( z)[\cdot])\,z - \nabla_\rho \mathcal
                       B(\rho( z)\nabla_ z \rho( z)[\cdot])\,u \bigr)\\
                  \nabla_\rho \mathcal C(\rho( z)\nabla_ z \rho( z)[\cdot])    % -\mathcal{C if C ist linear}
 \end{pmatrix} ,
\end{split} 
\end{equation}
where $G(z)[\cdot]$ stands for the linear map $G(z)\colon \mathbb R^{\ell}\to
\mathbb R^{\ell}$ so that e.g., $\lambda^\top G(z)[\cdot]$ can be interpreted as $(G(z)^\top[\lambda])^\top$. 

\begin{figure}[t]
         \centering
      \includegraphics[width=0.8\textwidth, height=0.43\textwidth]{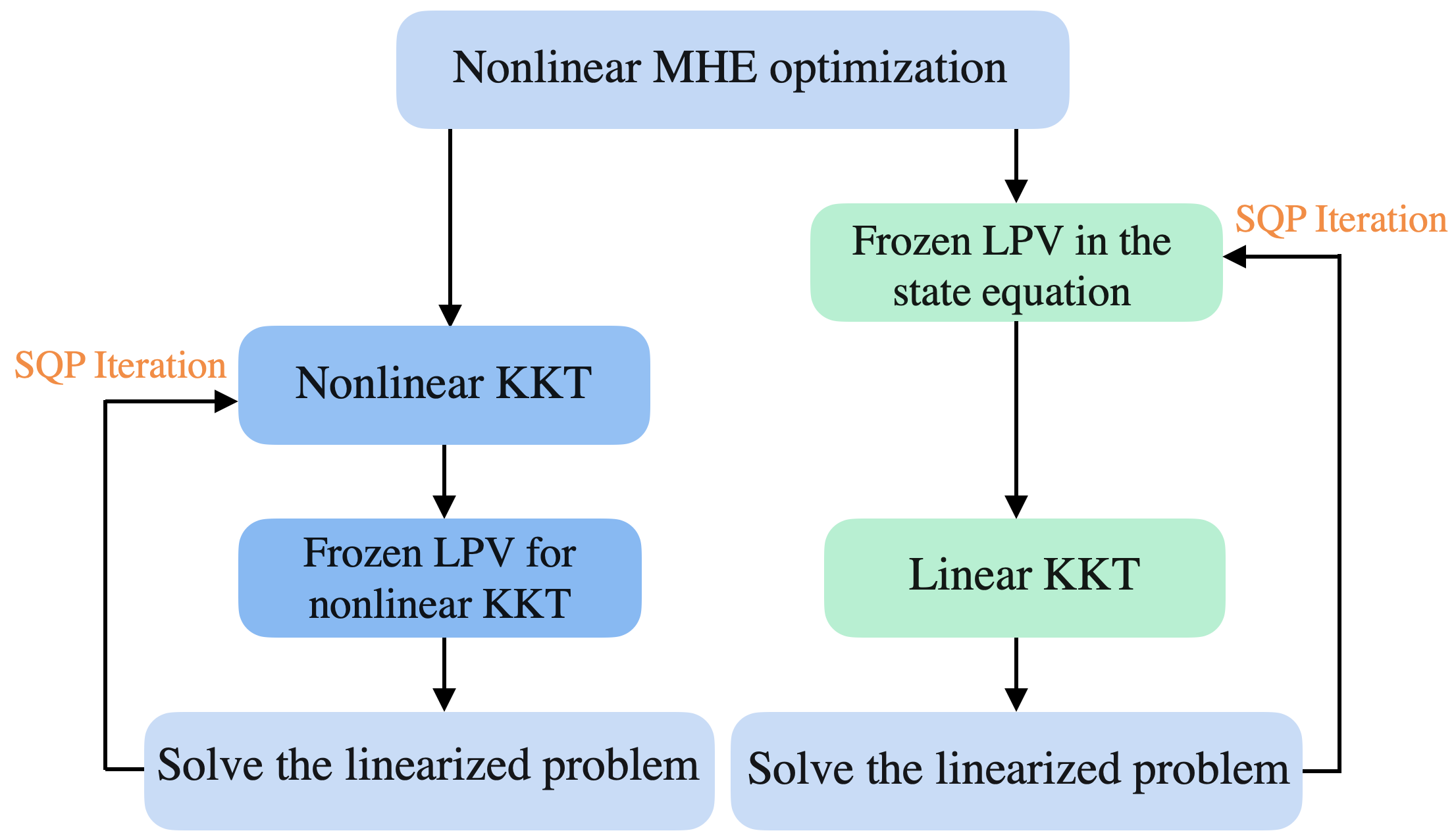} 
        \renewcommand{\figurename}{Figure}
      \caption{LPV-SQP-MHE scheme}
        \label{fig:two-SQP-MHE}
\end{figure}

Generally, the necessary optimality conditions (the KKT conditions) read as follows
\begin{subequations}\label{eq:kkt-full}
  \begin{align}
  h(z) &= 0, \label{eq:kkt-state} \\
  \lambda ^\top \nabla_z h(z) &= - \nabla_z f(z)\label{eq:kkt-costate}.
\end{align}
\end{subequations}
Note that \eqref{eq:kkt-state} resembles the state equations, while \eqref{eq:kkt-costate} defines the equations for the multiplier $\lambda$. These conditions form a coupled nonlinear system due to the state-dependent scheduling terms. A classical approach to solving such a system is to apply a  Newton-type method, which naturally leads to an SQP procedure. Each inner iteration, however, requires the Jacobians to be re-evaluated at the current iterate, which can be computationally expensive.
%solves a quadratic subproblem obtained from the linearization.
In this work, we exploit the quasi-LPV structure to avoid costly computations of Jacobian. By freezing the scheduling-dependent variables at the current iterate, the resulting subproblems become structured and directly lead to linear KKT conditions. 
% In particular, the computational cost per sampling time can be explicitly controlled by limiting the number of inner SQP iterations $M$.
% \begin{remark}
% In the numerical simulation, we select $M$ such that convergence error reach a stable and sufficiently small magnitude.
% \end{remark}
Based on where the scheduling parameter-dependent linearization is performed, we consider two approaches. The overview and comparison of both algorithms are given in Fig.\ref{fig:two-SQP-MHE}
%---------------------------------------Figure
%----------------------------------------
\subsection{Frozen LPV for state equation}
We start with the simplest approach.
At a given sampling time $k$, let $z^\ell$ denote the current SQP iterate, and we assume that a current trajectory $\rho^\ell =\rho(z^\ell)=\rho ( \{\bar x, \bar \omega, \bar \nu \})$ 
is provided to compute the updated state $z^{\ell+1}$ as the solution to the QP
\begin{equation*}
  \min_z f(z), \quad \text{s.t. } \mathcal A(\rho ^\ell)\,z - \mathcal
  B(\rho ^\ell))\,u=0, \quad y-\mathcal C(\rho ^\ell) z=0
\end{equation*}
and we define each QP subproblem as 
\begin{equation*}
        \min_{z}f(z)= \min_{z}
    \frac{1}2{}z^\top Hz+g^\top z ,
\end{equation*}
with the matrix $H$ and vector $g$ provided in the Appendix. Therefore, the constraints in \eqref{eq:kkt-full} become linear in the decision variables, i.e.
\begin{subequations}\label{eq:SQP-LPV-sys}
\begin{align}
  \mathcal A(\rho^\ell)\,z &=  \mathcal B(\rho^\ell) \,u , \\
  \mathcal C(\rho^\ell) z & =y, \\
 [ \mathcal A(\rho^\ell), \mathcal{C} (\rho^\ell)] \lambda  &=-( Hz+g). % \lambda ^T \mathcal A(\rho^\ell)  &= - \nabla_z f(z).
\end{align}
\end{subequations}
The obtained solution of $z$ defines the next iterate $z^{l+1}$, and this procedure is repeated until convergence or the maximum number of iterations $M$ is reached.

The above SQP iterations require an initial decision variable at each sampling instant. %At each sampling instant $k$, let $z_k^{(\ell)}$ denote the decision variable iterate at the \(\ell\)-th inner iteration.
For the first estimation window, $k=N$ set the initial guess $z^{0}$ to zero. For subsequent estimation windows $k>N$, the solution obtained in the previous sampling instant is used to initialize $z^{0}$. Specifically, the state estimates over the overlapping part of the horizon are shifted according to
\begin{equation*}
      \chi_{i|k}^{0}
    =
    \hat{x}_{k-N+i+1|k-1},
    \qquad i=0,\ldots,N-1,
\end{equation*}
where $\hat{x}_{k-N+i+1\,|\,k-1}$ denotes the estimated state at time
$k-N+i+1$ obtained from the previous MHE problem. The newly introduced terminal state is initialized by the prediction,
\begin{equation*}
      \chi_{N|k}^{0}
    =
    A\!\left(\rho(\chi_{N-1|k}^{0})\right)
    \chi_{N-1|k}^{0}
    +
    B\!\left(\rho(\chi_{N-1|k}^{0})\right)u_{k-1}.
\end{equation*}
The process noise and measurement noise estimates are initialized by shifting the corresponding solution from the previous horizon
\begin{equation*}
\begin{aligned}
    \omega_{i|k}^{0}
=
\hat{w}_{k-N+i+1|k-1},
\qquad i=0,\ldots,N-2\\
        \nu_{i|k}^{0}
=
\hat{v}_{k-N+i+1|k-1},
\qquad i=0,\ldots,N-1,
\end{aligned}
\end{equation*}
where $\hat{w}_{k-N+i+1\,|\,k-1}$ and $\hat{v}_{k-N+i+1|k-1}$ denote the estimated process noise and measurement noise at time $k-N+i+1$ obtained from the previous MHE problem. And the newly introduced noise values are initialized as
\begin{equation*}
    \begin{aligned}
       &  \omega_{N-1|k}^{0}=0\\
        &  \nu_{N|k}^{0}
    =
    y_k-
    C\!\left(\rho(\chi_{N|k}^{0})\right)
    \chi_{N|k}^{0}
    \end{aligned}
\end{equation*}
\begin{remark}
    Freezing the scheduling variables computed from the state trajectory of the current SQP iterate introduces an approximation of the original nonlinear dynamics within each SQP iteration. However, the approximation and the LPV model is refined throughout the optimization process, since the scheduling variables are recomputed after every iteration using the updated state trajectory. The present work does not provide theoretical bounds on the estimation error introduced by this scheduling update strategy. Instead, its effect is assessed numerically in Section \ref{sec:convergance test}. The numerical results show that the estimation error decreases as the number of SQP iterations increases, indicating that repeated scheduling updates effectively reduce the approximation error for the example considered.
\end{remark}
\begin{remark}
Although the resulting optimization problem can be interpreted as a QP problem, only equality constraints are considered in the formulation. Therefore, no general QP solver is required.
\end{remark}
%----------------------------------------
\subsection{Frozen LPV for nonlinear KKT}
\label{sec:Frozen LPV for nonlinear KKT}
Alternatively, we start from the nonlinear KKT conditions in \eqref{eq:kkt-full} and linearize them by freezing the nonlinear Jacobian in constraints at the current iterate. As the cost function was considered quadratic, the Jacobian in the constraints \eqref{eq:kkt-costate} is already linear for quadratic MHE, and the linearization will affect the term $\lambda^\top\nabla_z h(z)$. A frozen LPV formulation would read:
    \begin{equation*}
        \begin{aligned}
              \mathcal A(\rho^\ell)\,z &=  \mathcal B(\rho^\ell)\,u , \\
              \mathcal C (\rho^\ell) z & =y, \\
   \nabla_z h(z^\ell) ^\top\lambda   &= - (Hz+g). \label{eq:SQP-LPV-KKT-costate}
        \end{aligned}
    \end{equation*}
Compared to the first algorithm, this approach requires calculating the Jacobian of the constraints, and is expected to provide a more accurate approximation of the original nonlinear optimization conditions. However, due to the frozen LPV representation, the nonlinear terms are partially simplified before the Jacobian computation, leading to a structured and less complex approximation than in a conventional SQP framework.

%------------------------------------------------------------------------------------------------
\subsection{Newton-linearized MHE baseline}
\label{sec:Newton MHE}
To provide a linearization-based reference with a comparable iteration budget, we additionally consider a finite-iteration Newton-linearized MHE scheme.%, which is motivated by finite-iteration linearization-based MHE approaches \cite{KD11,AleG17}. 
In contrast to the proposed methods, it does not depend on scheduling-dependent LPV structure to construct the optimization update. This comparison therefore allows us to assess whether the computationally efficient estimation performance of the proposed methods results merely from limiting the number of optimization iterations or from exploiting the LPV structure in the update.

Starting from the KKT conditions in \eqref{eq:kkt-full}, we construct the Newton-linearized MHE update by applying a first-order linearization to the nonlinear equality constraints in the current MHE iterate. The Newton
linearization of the nonlinear constraints around $z^\ell$ is given by
\begin{equation*}
\label{eq:newton-constraint}
    h(z^\ell) +\nabla_z h(z^\ell)(z-z^\ell) \approx 0.
\end{equation*}
The KKT conditions of the corresponding quadratic subproblem are
\begin{subequations}
\label{eq:newton-KKT}
\begin{align}
    \nabla_z h(z^\ell)(z-z^\ell)
        &= -h(z^\ell),\label{eq:newton-KKT 1}
\\ \label{eq:newton-KKT 2}
        \nabla_z h(z^\ell) ^\top\lambda   &= - (Hz+g).
\end{align}
\end{subequations}
Note that the Jacobian is evaluated at the current iterate $z^l$  in the second equation  \eqref{eq:newton-KKT 2} and its value is fixed when solving for $z$ and $\lambda$. This equation is therefore identical to the second KKT equation in the LPV-based method in Section \ref{sec:Frozen LPV for nonlinear KKT}. Moreover, the Jacobian $\nabla_z h(z^\ell)$ required for the Newton linearization in the first equation \eqref{eq:newton-KKT 1} is already required by the second KKT equation. Therefore, the Newton-linearized method and the LPV-based method in Section \ref{sec:Frozen LPV for nonlinear KKT} have comparable computational costs in terms of the KKT condition, making the Newton-linearized method a suitable baseline for comparison under a similar computational budget.

%-------------------------------------------------------------------------------------------------------------------------------------
\section{Numerical Experiments}
\label{sec:simulation}
In this sequel, we will present the numerical nonlinear state estimation results for the two-link arm robot manipulator for the two presented LPV-SQP-MHE approaches. 
\subsection{The two-link arm robot model}
%----------------------------------------
We take the same two-link arm robot manipulator as in \cite{JiJM25}, which can be described using the following equations of motion 
\begin{equation}\label{2arm-nonlin model}
M(q(t))\ddot{q}(t) + c(q(t), \dot{q}(t))\dot{q}(t) + g(q(t)) = p_7\,\tau(t),
\end{equation}
where $q(t) := (q_1(t), q_2(t))$ are the joint angular positions and
$\tau(t):=(\tau_1(t),\tau_2(t))$ represents the motor torques. The model functions $M$, $c$, and $g$ are defined as
\begin{equation*}
  \begin{split}
    M(q(t)) &= \begin{bmatrix}
    \theta_1 & \theta_2 \cos(q_\Delta (t)) \\ 
    \theta_2 \cos(q_\Delta (t)) & \theta_3
\end{bmatrix}, \\
      g(q(t)) &= \begin{bmatrix}
    -\theta_4 \sin(q_1(t)) \\ 
    -\theta_5 \sin(q_2(t))
\end{bmatrix}, \\
        c(q(t), \dot q(t)) &= \begin{bmatrix}
    \theta_2 \sin(q_\Delta (t)) \dot{q}_2(t)^2 + \theta_6   \dot{q}_1(t) \\ 
    -\theta_2 \sin(q_\Delta(t)) \dot{q}_1(t)^2 + \theta_6 (\dot{q}_2(t) - \dot{q}_1(t))
\end{bmatrix}.
        \end{split}
\end{equation*}
with $q_\Delta(t) := q_1(t) - q_2(t)$.
The values of the parameter vector $\theta$ of the model are given in
Tab.~\ref{tab:params-twoarm-robot}. 
%----------------------------------
\begin{table}[t]
\centering
\begin{tabular}{|c|c|c|c|c|c|c|c|}
\hline
\textbf{Parameter:} & $\theta_1$ & $\theta_2$ & $\theta_3$ & $\theta_4$ & $\theta_5$ & $\theta_6$ & $\theta_7$ \\ \hline
\textbf{Value: (in the simulation)} & 5.6794 & 1.473 & 1.7985 & 0.4 & 0.4 & 2 & 1 \\ \hline
\textbf{Physical interpretation:} & 
\multicolumn{3}{|c|}{{Moments of inertia}} & 
\multicolumn{2}{|c|}{{Gravity scaling}} & 
{Friction} & 
{Input scaling} \\ \hline
\end{tabular}
\caption{Model parameters of the two-link robot arm.}
\label{tab:params-twoarm-robot}
\end{table}
%----------------------------------
Taking the torques as controls, we rewrite the model as a nonlinear
first-order state-space equation
\begin{equation}\label{eq:2arm-nonlin}
\begin{split}
  \dot{x}(t) &= \bar f(x(t)) + \bar B(x(t))\, u(t) \\
    \quad y(t) &= 
    \begin{bmatrix}
        I_{2 \times 2} & O
    \end{bmatrix} x(t),
    \end{split}
\end{equation}
with
\begin{equation*}\label{eq:2arm-nonlin-coefficients}
  \begin{split}
    \bar f(x) &= \begin{bmatrix}
      \begin{bmatrix}
    x_3(t) \\ x_4(t)
\end{bmatrix} \\
        M\left(\begin{bmatrix}
    x_1(t) \\ x_2(t)
\end{bmatrix} \right )^{-1}g \left (\begin{bmatrix}
    x_1(t) \\ x_2(t)
\end{bmatrix} \right ) - M \left (\begin{bmatrix}
    x_1(t) \\ x_2(t)
\end{bmatrix} \right )^{-1}c \left( \begin{bmatrix}
    x_1(t) \\ x_2(t)
\end{bmatrix},  \begin{bmatrix}
    x_3(t) \\ x_4(t)
\end{bmatrix} \right ) \,
        \end{bmatrix},\\
      \bar B(x) &= \begin{bmatrix}
        0 \\ p_7 M \left(\begin{bmatrix}
    x_1(t) \\ x_2(t)
\end{bmatrix} \right )^{-1}
      \end{bmatrix} \\
  \end{split}
\end{equation*}
and where $x(t) = \left (\begin{bmatrix}
    x_1(t) \\ x_2(t)
\end{bmatrix}, \begin{bmatrix}
    x_3(t) \\ x_4(t)
\end{bmatrix} \right ) ^\top := (q(t), \dot{q}(t)) ^\top$ and $u(t) =
\tau(t)$. In the LPV notation that follows, we will consider $x(t)$ as a column
vector, i.e.
\begin{equation}\label{eq:xascvec}
  x(t) = \begin{bmatrix}
   x_1(t)\\x_2(t)\\x_3(t)\\x_4(t)
\end{bmatrix}=\begin{bmatrix}
   q_1(t)\\ {q}_2(t)\\\dot{q}_1(t)\\\dot{q}_2(t)
\end{bmatrix}.
\end{equation}

\[
x(t)
=
\begin{bmatrix}
\begin{bmatrix}
    x_1(t) \\ x_2(t)
\end{bmatrix} \\
\begin{bmatrix}
    x_3(t) \\ x_4(t)
\end{bmatrix} \\ 
\end{bmatrix}
:=
\begin{bmatrix}
q(t)\\
\dot q(t)
\end{bmatrix},
\qquad
u(t)=\tau(t).
\]
%-----------------------switched to LPV form
\subsection{The Model in Quasi LPV representation}
We consider a so-called \emph{scheduling map} $\rho(t)=\rho(x(t))=\rho((q(t),\dot q(t)))$ chosen as in \cite{KoeT20} as $\rho: \mathbb{R}^4 \to \mathbb{R}^{10}$:
\begin{equation}\label{eq:lpv-rho-sched-map}
  % \begin{split}
    \rho((q,\dot{q})) = % \\
            \begin{bmatrix}
        \frac 1h \\
        \frac 1h \cos(q_\Delta) \\
        \frac 1h \sin(q_1)\frac{1}{q_1} \\
        \frac 1h \cos(q \Delta)\sin(q_2)\frac{1}{q_2} \\
        \frac 1h (- \theta_2^2 \sin(q_\Delta) \cos (q_\Delta) \dot q_1 - \theta_6(\theta_3 +
        \theta_2\cos(q_\Delta))) \\
        \frac 1h (-\theta_3\sin(q_\Delta) \dot{q_2} +  \theta_6\cos(q_\Delta)) \\
        \frac 1h \cos(q_\Delta) \sin(q_1)\frac{1}{q_1} \\
        \frac 1h \sin(q_2)\frac{1}{q_2} \\
        \frac 1h (\theta_1\theta_2\sin(q_\Delta)\dot{q_1} + \theta_6(\theta_1 + \theta_2\cos(q_\Delta))) \\
      \frac 1h (\theta_2^2\sin(q_\Delta)\cos(q_\Delta) \dot{q}_2 - \theta_1\theta_6)
    \end{bmatrix},
  % \end{split}
\end{equation}
where $h = \operatorname{det}(M) = \theta_1\theta_3 - \theta_2^2\cos^2(q_\Delta)$.
We note that the function $\frac{\sin(s)}{s}$ is continuously extendable at
$s=0$ and that with the chosen parameters (see Tab.
\ref{tab:params-twoarm-robot}) the value of $1/h$ is bounded.
With the help of \eqref{eq:lpv-rho-sched-map}, the nonlinear model \eqref{eq:2arm-nonlin} with $x$ as defined in
\eqref{eq:xascvec} can be equivalently rewritten as a so-called quasi LPV system:
\begin{equation}\label{eq:2arm-nonlin LPV form}
\begin{aligned}
    \dot{x}(t) &= A(\rho(t))x(t) + B(\rho(t))u(t), \\
    y(t) &= Cx(t) + Du(t),
    \end{aligned}
\end{equation}
with
\begin{equation*}
  \begin{split}
    A(\rho) &= 
    \begin{bmatrix}
        0 & 0 & 1 & 0 \\
        0 & 0 & 0 & 1 \\
       - \theta_3 \, \theta_4\,\rho_3 & \theta_2 \, \theta_5\,\rho_4 & \rho_5 & \theta_2\,\rho_6 \\
        \theta_2 \, \theta_4\,\rho_7 & - \theta_1\, \theta_5\,\rho_8 & \rho_9 & \rho_{10}
    \end{bmatrix},\\
    B(\rho) &= 
    \begin{bmatrix}
        0 & 0 \\
        0 & 0 \\
        \theta_3\, \theta_7\,\rho_1 & -\theta_2\, \theta_7\,\rho_2 \\
        -\theta_2\, \theta_7\,\rho_2 & \theta_1\, \theta_7\,\rho_1
    \end{bmatrix}, \quad
    C = \begin{bmatrix}
      I_{2 \times 2} & O
    \end{bmatrix}, \quad D = 0.
  \end{split}
\end{equation*}
Here, all nonlinearities are factorized into the coefficient matrices $A$ and $B$.
%-----------------------------------------------
\subsection{Comparison of two LPV-SQP-MHE approaches}
We now evaluate the proposed LPV-SQP-MHE schemes on the two-link arm robot manipulator. The continuous-time system is discretized using the forward Euler method with sampling time $dt=0.1s$, and the MHE is performed over the total simulation time $T=20s$ ($N_f=\frac{T}{dt}=200$ steps) with a fixed sliding window of $T_n=1s$ ($N=\frac{T_n}{dt}=10$ steps). The input signal is selected as $u(t)=[\sin(t),\cos(t)]^\top $. We begin the state estimation from the time step $k=N$. 
The two proposed LPV-SQP-MHE schemes are compared with the Newton-linearized MHE introduced in Section \ref{sec:Newton MHE} and with a nonlinear MHE solved with CasADi/IPOPT. The Newton-linearized method provides a finite-iteration linearization-based baseline with a computational budget comparable to that of the proposed methods. Performance is assessed in terms of state estimation error and computational runtime for the different number of SQP iterations $M$. 
In the simulations, the true initial state is \(x(0)=[0.1,0.1,0.1,0.1]^\top\), whereas the observers are initialized at $\bar x(0)=[0,0,0,0]^\top$. After the first estimation step, the solution from the previous time step is used to initialize the next MHE problem, which is applied to all methods. The process and measurement noises are generated as zero-mean Gaussian sequences. Specifically, we draw $w_k\sim\mathcal N(0,0.1 \cdot I_{4 \times 4})$ and $v_k\sim\mathcal N(0,0.1 \cdot I_{2 \times 2})$.

\begin{remark}
MHE does not require Gaussian noise assumption, the disturbances can be modeled flexibly, for example, as bounded noise. However, Kalman-filter based methods, require zero-mean Gaussian assumption. For the potential future comparison, we therefore generate noise as zero-mean Gaussian in the simulation.
\end{remark}
% -------------------new command %use \texttt{runtime} directly...
%\newcommand\runtime{\texttt{CT}}
%\runtime 

To quantify the simulation of estimation accuracy, we compute the state-trajectory error
\begin{equation*}
e_x=\|X_{\text{est}}-X_{\text{real}}\|_F,
\qquad
X_{\text{est}}=[\hat x_0,\hat x_1,\dots,\hat x_{N_f}],
\end{equation*}
where \(X_{\mathrm{est}}\) and \(X_{\mathrm{real}}\) denote, respectively, the estimated state trajectory and the true state trajectory generated by the system dynamics with the process and measurement noise over the total simulation time $T$. 
$\texttt{runtime}$ denotes the average running time over five simulations in the interval $[0, T]$. We denote the scheme in which we directly apply the frozen LPV parameters to the state equations by \texttt{fLPV\_state}, and the scheme based on freezing the LPV parameters in the nonlinear KKTs by \texttt{fLPV\_KKT}.

Now we compare the proposed LPV-SQP-MHE schemes with the Newton-linearized MHE and the nonlinear MHE solved by CasADi/IPOPT. The estimation is then updated at each time step with a fixed horizon length $N$.
The state estimation error and the runtime under the sampling time $dt=0.1s$ are summarized in Tab.~ \ref{tab:error and runtime}. The values are obtained by averaging over five independent runs with different sets of noise. For the LPV-SQP-MHE methods and Newton-linearized MHE, we perform one and five iterations, i.e. $M=1$ and $M=5$, whereas the CasADi/IPOPT implementation is solved to convergence at each sampling step. Note that small $dt$ increases the number of estimation updates $N_f=T/dt$ over a fixed duration $T$, which leads to a higher total runtime.
For the proposed LPV-SQP-MHE schemes, the estimation error remains relatively small even with $M=1$, and the estimation errors of the two schemes are almost identical. For both $M=1$ and $M=5$, the proposed methods and the Newton-linearized MHE require very similar computational $\texttt{runtime}$, while the proposed methods achieve smaller estimation errors. % The results indicate that, for the considered example, the LPV-based method achieves higher estimation accuracy than the Newton linearization under a comparable computational budget.
Compared to the nonlinear MHE solved using CasADi, the estimation error is only one order of magnitude larger. In terms of computational cost, the proposed methods achieve comparable estimation performance with significantly reduced computational runtime. %the proposed methods are significantly faster than the nonlinear MHE solved by CasADi/IPOPT.

%-----------------------------------table M=1

\begin{table}[t]
\begin{center}
\renewcommand{\arraystretch}{1.1}
\resizebox{\textwidth}{!}{
\begin{tabular}{|c|c|c|c|c|c|c|c|}
\hline
& \multicolumn{2}{|c|}{\textbf{ \texttt{fLPV\_state}} } & \multicolumn{2}{|c|}{\textbf{ \texttt{fLPV\_KKT}} }     & 
\multicolumn{2}{|c|}{\textbf{ \texttt{Newton-linerized MHE}} } 
& \textbf{CasADi MHE}   \\ \hline
& $M=1$&$M=5$&$M=1$&$M=5$ & $M=1$ & $M=5$ & Full iteration  \\ \hline
\textbf{$e_x$} &   \begin{tabular}{c}
       $2.2214\times 10^{-03}$
   \end{tabular}  & \begin{tabular}{c}
       $1.1217\times 10^{-03}$
   \end{tabular}  &  \begin{tabular}{c}
       $2.2214\times 10^{-03}$
   \end{tabular}&\begin{tabular}{c}
       $1.1215\times 10^{-03}$
   \end{tabular} &   \begin{tabular}{c}
       $1.6125\times 10^{-02}$
   \end{tabular}  & \begin{tabular}{c}
       $1.2037\times 10^{-02}$
   \end{tabular}  
   &$7.2521\times 10^{-04}$ \\ \hline
% NLPV observer of~\cite{LPV_2020_RNC} & \multicolumn{8}{|c|}{Infeasible solution} \\ \hline
% ~\textcolor{red}{LPV Observer}~%\cite{} & \multicolumn{8}{|c|}{Infeasible solution} \\ \hline
$\textbf{\texttt{runtime}} $~%\cite{} 
& $0.3014 s$& $1.2192 s$& $1.8692 s$&$2.0524 s$ & $1.9062 s$  & $2.0692 s$  & $28.1725 s$   \\ \hline

\end{tabular}
}
%\vspace{-0.5cm}
\end{center}
\caption{$\mathcal L_2$ values of estimation error $e_x$ and \texttt{runtime}
for the different MHE schemes on sampling tme $dt=0.1s$.}
\label{tab:error and runtime}
\end{table}

\begin{figure}[t]
         \centering
      \includegraphics[width=0.83\textwidth, height=0.48\textwidth]{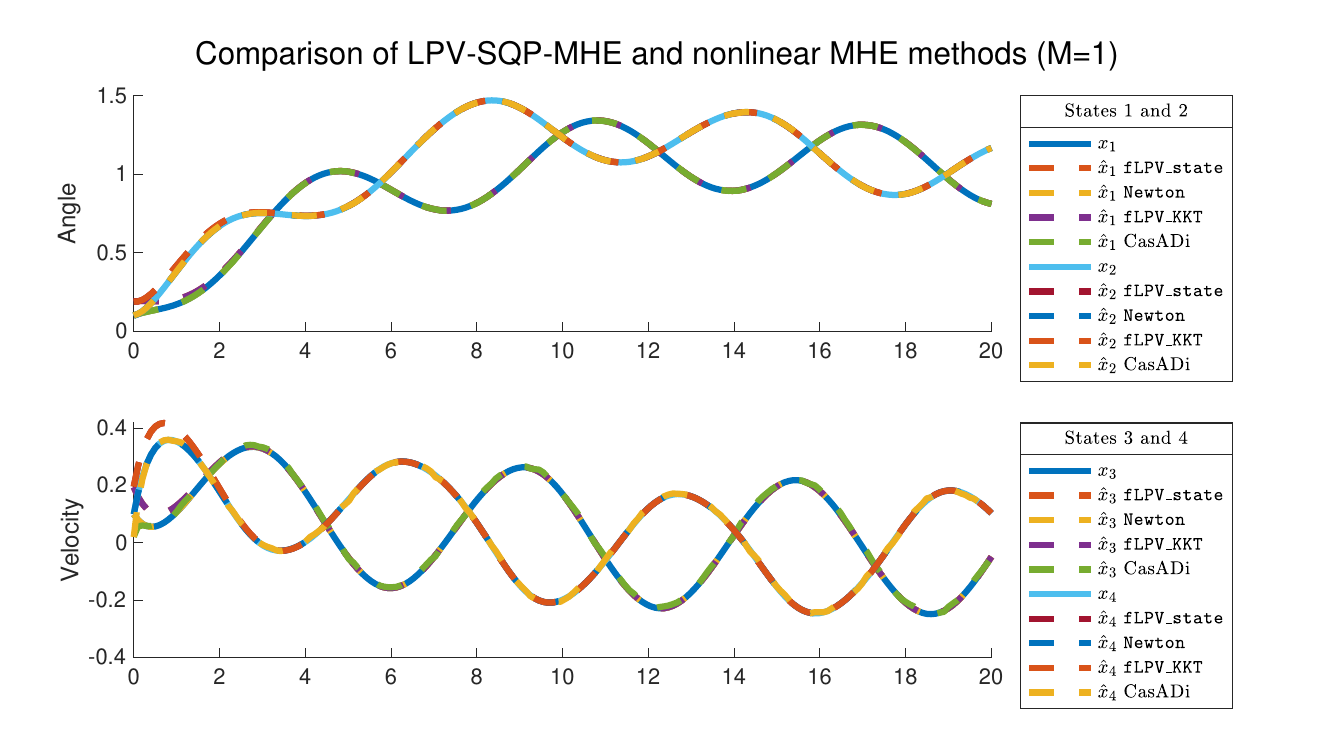} 
        \renewcommand{\figurename}{Figure}
      \caption{Comparison of state estimation trajectories for different MHE schemes.}
        \label{fig:all comparison}
\end{figure}

%---------------------------------------Figure
The resulting estimated state trajectories obtained in one of the simulations are shown in Fig.~ \ref{fig:all comparison}, corresponding to a sampling time of $dt=0.1s$ and one SQP iteration $M=1$. % For the LPV-SQP-MHR methods, we perform one and five iterations, i.e., $M=1$ and $M=5$. 
All methods can effectively reproduce the trajectory of the true state. For the nonlinear MHE solved by CasADi/IPOPT, even if the initial prior estimate of the state deviates from the true value, the solver can correct this error within a single estimation step, since the internal optimization is fully iterated at each step and nonlinear KKTs in the optimization problem are not linearized. 
%Compared to the nonlinear MHE implementation solved using CasADi/IPOPT, the proposed LPV-SQP-MHE methods only perform a limited number of SQP iterations. Therefore, small deviations may appear at the beginning of the curve. 
Compared to the nonlinear MHE implementation solved using CasADi/IPOPT, which solves the nonlinear optimization problem using full nonlinear Jacobian and Hessian, the proposed LPV-SQP-MHE methods rely on LPV-based approximation and maximum number of iterations $M$. Therefore, small deviations may appear at the beginning of the estimated trajectories.
However, as new measurement data are incorporated over time, the estimates are continuously updated toward the true trajectory. As shown in the figure, the estimated trajectory quickly aligns with the true state trajectory from the beginning. A similar behavior can also be observed for the Newton-linearized MHE.
%-------------------------------------------------------------------------
\subsection{Convergence tests}
\label{sec:convergance test}
To numerically analyze the convergence behavior of the proposed schemes, we consider an estimation window at $k=N$ with $M=20$, and evaluate the same set of state trajectories under different SQP iterations. Since the initial prior is chosen different from the true initial state, it is clear to see how the SQP iterations reduce the initial estimation error in the first estimation window. First, we collect the estimated state trajectory within the first window as
\begin{equation*}
    X^\ell = [x_0^\ell,x_1^\ell,\dots,x_N^\ell]
\end{equation*}
and then define the inner SQP iteration error as
\begin{equation*}
    e^\ell = \|X^{\ell+1}-X^\ell\|_F .
\end{equation*}
%--------------------- convergence picture
%  \begin{figure}[h]
%          \centering
%       \includegraphics[width=0.91\textwidth, height=0.32\textwidth]{image/LPV_SQP_MHE_convergence.pdf} 
%         \renewcommand{\figurename}{Figure}
%       \caption{Convergence of the SQP iteration error for a horizon window $k=N$}
%         \label{fig:convergence figure}
% \end{figure}
%--------------------convergence table
\begin{table}[t]
\centering
\begin{tabular}{|c|c|c|c|c|c|c|}
\hline
\textbf{\textbf{$M$}} & $1$ & $2$ & $5$ & $10$ & $20$  \\ \hline % & $6$ & $7$ &$8$ &$9$ &$10$ 
\textbf{ \texttt{fLPV\_state}}  &$ 1.1424$ & $8.1674\times 10^{-06}$ & $6.4260\times 10^{-06}$ & $4.1830\times 10^{-06}$& $1.7881\times 10^{-06}$ \\ \hline % \texttt{fLPV_{state}}
\textbf{\texttt{fLPV\_KKT}}  & $1.1424$& $8.1066\times 10^{-06}$ & $6.3915\times 10^{-06}$ & $4.1750\times 10^{-06}$& $1.7970\times 10^{-06}$ \\ \hline % \texttt{fLPV_{KKT}}
\end{tabular}
\caption{Convergence of the SQP iteration error}
\label{tab: Convergence table}
\end{table}

Similarly, we collect the mean values of the iteration error from five different sets of noise. The detailed numerical values of the iteration error are summarized in Tab. ~\ref{tab: Convergence table} (For readability, we only show the results of SQP iteration numbers \(M=1,2,5,10,20\)). As the number of iterations $M$ increases, the iteration error consistently decreases for both proposed methods, showing numerically that the iterates approach a stable solution for the considered case. A significant reduction in the error is already observed after the first iteration, and from the second iteration, only small changes between successive iterates are observed.

% Due to the large difference in magnitude between the first iteration error and the others, the Fig. ~\ref{fig:convergence figure} focuses on the number of iteration $M=2$ to $M=20$ to better illustrate the convergence trend.  

%-------------------------------------------------------------------------------------------------------------------------------------
\section{Conclusion}
In this work, we have designed two LPV-SQP-MHE schemes based on different linearization strategies, and applied them to a two-link robot model. Both methods reformulate the nonlinear MHE problem into structured quadratic subproblems through LPV-based linearization. This allows the Jacobian computations to be replaced by LPV updates.
The two proposed LPV-SQP formulations are shown to reduce the computational burden by replacing the nonlinear MHE subproblem with a sequence of  structured KKTs. Compared to the nonlinear MHE solved by CasADi/IPOPT, they significantly  reduce computational costs while still maintaining estimation performance. Moreover, the comparison with the finite-iteration Newton-linearized MHE shows that, for the considered example and under a comparable iteration budget, the proposed LPV-based schemes achieve smaller estimation errors.

As a future work, comparing the proposed method with estimation techniques such as Kalman-filter-based methods by incorporating more complex noise models would facilitate a more comprehensive evaluation. Furthermore, the current framework can be extended to more complex, higher-dimensional nonlinear systems, and to provide larger computational benefits.

\label{sec:conclusion}

\label{sec:conclusion}
%-------------------------------------------------------------------------------------------------------------------------------------
\section{Appendix}
According to the constraints on the state and inputs, the matrices $\mathcal A(\rho), \mathcal B(\rho), \mathcal C(\rho)$ are defined by 
\begin{equation*}
    \mathcal A(\rho)=\begin{bmatrix} \begin{array}{ccccc|cccc|ccc}
 -A(\rho(\chi_{k-N})) &  I & 0 & \cdots & 0 &  I &  0 & \cdots &   0 &    0 &\cdots &0 \\
        0 &   -A(\rho(\chi_{k-N+1})) &  I & 0  & \cdots     &  I & 0 & \cdots &  0 &0&\cdots &0\\
        \vdots    &  \ddots  &  \ddots  &  \ddots  &  \ddots  &  \vdots  &  \ddots  &  \ddots  &  \ddots  &  \vdots  &  \cdots &  \vdots
    \end{array}
    \end{bmatrix},
\end{equation*}
\begin{equation*}
    \mathcal B (\rho)= \begin{bmatrix}
         B(\rho(\chi_{k-N})) \\ B(\rho(\chi_{k-N+1})) \\ \vdots \\ B(\rho(\chi_{k-1}))
    \end{bmatrix}, 
    \quad \mathcal C (\rho)= \begin{bmatrix} \begin{array}{cccc|ccc|cccc}
        C(\rho(\chi_{k-N})) &  0 &  0 & \cdots    &  0  & \cdots & 0 &     I  & 0 & \cdots  & \cdots  \\
          0 & C(\rho(\chi_{k-N+1})) & 0 & \cdots    & \vdots & \cdots & \vdots &             0 &   I &   0 & \cdots                             \\
         \vdots & \ddots & \ddots & \ddots    & \vdots & \cdots & \vdots &              \vdots & \ddots & \ddots & \ddots
    \end{array}
    \end{bmatrix}.
\end{equation*}
\vspace{0.1 cm}
For a general horizon length $N$, we consider the decision variables $z = [x_0^\top,\dots,x_N^\top,\omega_0^\top,\dots,\omega_{N-1}^\top,
\nu_0^\top,\dots,\nu_N^\top]^\top.$ The MHE cost function $f(z)$ can be expressed in the quadratic form
\begin{equation*}
    f(z)=\frac12 z^\top H z + g^\top z + \text{const},
\end{equation*}
where
\begin{equation*}
    H = \mathrm{blkdiag}(
P, 0,\dots,0,Q,\dots,Q,R,\dots,R ), \quad 
g =
\begin{bmatrix}
- P\bar{x}_0 \\
0\\
\vdots \\
0
\end{bmatrix},
\end{equation*}
where $\mathrm {blkdiag}$ denotes a block diagonal matrix.

%--------------------------------------------------
% \begin{thebibliography}{99}

% \bibitem{1} Spiegel, M. R. (1981). Theory and problems of Advanced Calculus: Si (metric) edition. McGraw-Hill. 

% \end{thebibliography}
\bibliographystyle{abbrv}
\bibliography{LPV_SQP_MHE.bib}

\end{document}